\documentclass[12pt,a4paper]{article}
\usepackage{amsfonts}
\usepackage{amssymb}
\usepackage{amsmath}

\newtheorem{theorem}{Theorem}[section]
\newtheorem{definition}[theorem]{Definition}

\newtheorem{corollary}[theorem]{Corollary}

\date{}

\begin{document}

\title{Gr\"{o}bner-Shirshov bases for
some Lie algebras\footnote{Supported by the NNSF of China
(11171118), the Research Fund for the Doctoral Program of Higher
Education of China (20114407110007), the NSF of Guangdong Province
(S2011010003374) and the Program on International Cooperation and
Innovation, Department of Education, Guangdong Province
(2012gjhz0007).}}
\author{
 Yuqun
Chen, Yu Li  and Qingyan Tang
\\
\\
{\small \ School of Mathematical Sciences, South China Normal
University}\\
{\small Guangzhou 510631, P. R. China}\\
{\small  yqchen@scnu.edu.cn}\\
{\small  LiYu820615@126.com}\\
{\small tangqingyan0812@163.com}}

\maketitle

\noindent\textbf{Abstract:} We give Gr\"{o}bner-Shirshov bases for
Drinfeld-Kohno Lie algebra  $\textbf{L}_{n}$ in \cite{[Et]} and
Kukin Lie algebra $A_P$ in \cite{Kukin}, where $P$ is a semigroup.
As applications, we show that as $\mathbb{Z}$-module
$\textbf{L}_{n}$ is free and a $\mathbb{Z}$-basis of
$\textbf{L}_{n}$ is given. We give another proof of Kukin Theorem:
if semigroup $P$  has the undecidable word problem then the Lie
algebra $A_P$ has the same property.



\section{Introduction}

Gr\"{o}bner bases and Gr\"{o}bner-Shirshov bases  were invented
independently by A.I. Shirshov for ideals of free (commutative,
anti-commutative) non-associative algebras \cite{Sh62a,Shir3}, free
Lie algebras \cite{Sh62b,Shir3} and implicitly free associative
algebras \cite{Sh62b,Shir3}  (see also \cite{Be78,Bo76}), by H.
Hironaka \cite{Hi64} for ideals of the power series algebras (both
formal and convergent), and by B. Buchberger \cite{Bu70} for ideals
of the polynomial algebras.

Gr\"{o}bner bases and Gr\"{o}bner-Shirshov bases theories have been
proved to be very useful in different branches of mathematics,
including commutative algebra and combinatorial algebra, see, for
example, the books \cite{AL, BKu94, BuCL, BuW, CLO, Ei}, the papers
\cite{Be78,Bo72,Bo76,BCC08,BCD08,BCLi13,BCL08,BCM,BCM13,
Bu70,CC12,KL,MZ},
 and the surveys \cite{BC,BC13,BCS, BFKK00, BK03, BK05}.

A.A. Markov \cite{Markov}, E. Post \cite{Post46}, A. Turing
\cite{Turing}, P.S. Novikov \cite{Novikov} and W.W. Boone
\cite{Boon} constructed finitely presented semigroups and groups
with the undecidable word problem. This result also follows from the
Higman  theorem \cite{Higman} that any recursive presented group is
embeddable into finitely presented group. A weak analogy of Higman
theorem for Lie algebras was proved in \cite{Bo72} that was enough
for existence of a finitely presented Lie algebra with the
undecidable word problem.

In \cite{Kukin}, Kukin constructed  the Lie algebra $A_P$ for a
semigroup $P$ such that if $P$ has the undecidable word problem then
$A_P$ has the same property. In this paper we  find a
Gr\"{o}bner-Shirshov basis for $A_P$ and then give another proof of
the above result.

The Drinfeld-Kohno Lie algebra $\textbf{L}_n$ appears in \cite{[Dr],
[Ko]} as the holonomy Lie algebra of the complement of the union of
the diagonals $z_{i}=z_{j},\ i<j$. The universal
Knizhnik-Zamolodchikov connection takes values in this Lie algebra.
In this paper, we give a Gr\"obner-Shirshov basis for the
Drinfeld-Kohno Lie algebra $\textbf{L}_n$ over  $\mathbb{Z}$ and a
$\mathbb{Z}$-basis of $\textbf{L}_n$. As an application we get a
simple proof that  $\textbf{L}_{n}$ is an iterated semidirect
product of free Lie algebras.

We are very grateful to Professor L.A. Bokut for  his  guidance and
useful discussions.

\section{Composition-Diamond lemma for Lie algebras over a field}

For the completeness of the paper, we formulate the
Composition-Diamond lemma for Lie algebras over a field in this
section, see \cite{BC13,Lyndon, Sh58,Shir3} for details.

Let $k$ be  a filed, $I$ a well ordered index set, $X=\{x_i |i\in
I\}$ a set, $X^*$ the free monoid generated by $X$ and $Lie(X)$ the
free Lie algebra over $k$ generated by $X$.

We order $X=\{x_i |i\in I\}$ by $x_i>x_t$ if $i>t$ for any $i,t\in
I$.

 We use two linear orders on $X^*$: for any $u,v\in X^*$,

(i) ({\sl lex order}) $1\succ t$ if $t\neq1$ and, by induction, if
$u=x_{i}u_{i}$ and $v=x_{j}v_{j}$ then $u\succ v$ if and only if
$x_{i}>x_{j}$, or $x_{i}=x_{j}$ and $u_{i}\succ v_{j}$;

(ii) ({\sl deg-lex order}) $u >v$ if and only if $deg(u)>deg(v)$, or
$deg(u)=deg(v)$ and $u\succ v$, where $deg(u)$ is the length of $u$.

We regard $Lie(X)$ as the Lie subalgebra of the free associative
algebra $k \langle X \rangle$, which is generated by $X$ under the
Lie bracket $[u,v]=uv-vu$. Given $f \in k \langle X \rangle$, denote
$\bar{f}$ the leading word of $f$ with respect to the deg-lex order.
$f$ is  {\sl monic} if the coefficient of $\bar{f}$ is 1.

\begin{definition}An associative
word $w\in X^*\setminus\{1\}$ is an associative Lyndon-Shirshov
 word (ALSW for short) if
$$
(\forall u,v\in X^*, u,v\neq1) \ w=uv\Rightarrow w>vu.
$$

A non-associative word $(u)$ in $X$ is a non-associative
Lyndon-Shirshov word (NLSW for short), denoted by $[u]$,  if

(i) $u$ is an ALSW;

(ii) if $[u]=[(u_{1})(u_{2})]$ then both $(u_{1})$ and $(u_{2})$ are
NLSW's;

(iii) if $[u]=[[[u_{11}][u_{12}]][u_{2}]]$ then $u_{12} \preceq
u_{2}$.
\end{definition}
We denote the set of all ALSW's and  NLSW's in $X$ by
 $ALSW(X)$ and $NLSW(X)$ respectively.

For an ALSW $w$, there is a unique  bracketing $[w]$ such that $[w]$
is NLSW:  $[w]=[[u][v]]$ if $deg(w)>1$, where $v$ is the longest
proper associative Lyndon-Shirshov end of $w$.

\ \

\noindent{\bf Shirshov Lemma} Suppose that $w=aub$, where $w,u\in
ALSW(X)$. Then
\begin{enumerate}
\item[(i)]
$ [w]=[a[uc]d], $ where  $b=cd$ and possibly $c=1$.
\item[(ii)] Represent $c$
in the form $ c=c_{1}c_{2} \ldots c_{n}, $ where $c_{1}, \ldots
,c_{n}\in ALSW(X)$ and $c_{1} \leq c_{2} \leq \ldots \leq c_{n}$.
Replacing $[uc]$ by $[\ldots[[u][c_{1}]] \ldots [c_{n}]]$ we obtain
the word
$$
[w]_{u}=[a[\ldots[[[u][c_{1}]][c_2]] \ldots [c_{n}]]d]
$$
which is called the {\it Shirshov special bracketing} of $w$
 relative to  $u$.
\item[(iii)] $ \overline{[w]}_{u}=w. $
\end{enumerate}

\begin{definition}
Let $S\subset Lie(X)$ with each $s\in S$ monic, $a,b\in{X^*}$ and
$s\in S$. If $a\bar{s}b$ is an ALSW, then we call
$[asb]_{\bar{s}}=[a\bar{s}b]_{\bar{s}}|_{[\bar{s}]\mapsto{s}}$ a
special normal $S$-word, where $[a\bar{s}b]_{\bar{s}}$ is defined as
in Shirshov Lemma. An $S$-word $(asb)$ is called a normal $S$-word
if $\overline{(asb)}=a\overline{s}b$.

Suppose that $f,\ g\in S$. Then, there are two kinds of
compositions:
\begin{enumerate}
\item[(i)] If  $w=\bar{f}=a\bar{g}b$ for some $a,b\in X^*$, then the
polynomial $( f,g)_w=f - [agb]_{\bar{g}}$ is called the
 composition of inclusion of $f$ and $g$ with respect to $w$.

\item[(ii)] If \ $w$ is a word such that $w=\bar{f}b=a\bar{g}$ for
some $a,b\in X^*$ with $deg(\bar{f})+deg(\bar{g})>deg(w)$, then the
polynomial
 $( f,g)_w=[fb]_{\bar{f}}-[ag]_{\bar{g}}$ is called the composition of intersection of $f$ and
$g$ with respect to $w$.
\end{enumerate}

In (i) and (ii), $w$ is called an ambiguity.

Let  $h$ be a Lie polynomial and $w\in X^*$. We shall say that $h$
is trivial modulo $(S,w)$, denoted by $h\equiv_{Lie}0 \  mod(S,w)$,
if $h=\sum_{i}\alpha_i(a_is_ib_i)$, where each $(a_is_ib_i)$ is a
normal $S$-word and $a_i\bar{s_i}b_i<w$.

The set $S$ is called a  Gr\"obner-Shirshov basis in $Lie(X)$ if any
composition
 in $S$ is trivial modulo $S$ and corresponding $w$.
\end{definition}

\begin{theorem}\label{cdLie} ({\bf Composition-Diamond lemma for Lie
algebras over a field}) Let $S\subset{Lie(X)}\subset{k\langle
X\rangle}$ be nonempty set of monic Lie polynomials. Let $Id(S)$ be
the ideal of $Lie(X)$ generated by $S$. Then the following
statements are equivalent.
\begin{enumerate}
\item[(i)] $S$ is a Gr\"{o}bner-Shirshov basis in
$Lie(X)$.
\item[(ii)] $f\in{Id(S)}\Rightarrow{\bar{f}=a\bar{s}b}$ for
some $s\in{S}$ and $a,b\in{X^*}$.
\item[(iii)]$Irr(S)=\{[u]\in NLSW(X) \ | \  u\neq{a\bar{s}b}, \
s\in{S},\ a,b\in{X^*}\}$ is a linear basis for
$Lie(X|S)=Lie(X)/Id(S)$.
\end{enumerate}
\end{theorem}

\noindent{\bf Remark:} The above Composition-Diamond lemma is also
valid if we replace the base field $k$ by an arbitrary commutative
ring $K$ with identity. If this is the case, as $K$-module,
$Lie(X|S)$ is free with a $K$-basis $Irr(S)$.

\section{Kukin's construction of a Lie algebra with unsolvable word problem}

Let $P=sgp\langle x,y|u_i=v_i, \ i\in I\rangle$ be a semigroup.
Consider the Lie algebra
$$
A_P=Lie(x, \hat{x}, y, \hat{y}, z|S),
$$
where $S$ consists of the following relations:
\begin{enumerate}
\item[(1)] $[\hat{x}x]=0,\ [ \hat{x}y]=0,\ [ \hat{y}x]=0,\ [\hat{y}y]=0$,
\item[(2)] $[\hat{x}z]=-[zx],\  [\hat{y}z]=-[zy]$,
\item[(3)] $\lfloor zu_i\rfloor=\lfloor zv_i\rfloor,\ i\in  I$.
\end{enumerate}
Here, $\lfloor zu\rfloor$ means the left normed bracketing.

In this section, we give a Gr\"{o}bner-Shirshov basis for Lie
algebra $A_P$ and by using this result we give another proof for
Kukin's theorem, see Corollary \ref{t4.12}.

\ \

Let the order $ \hat{x}>  \hat{y}> z>x>y$ and $>$   the deg-lex
order on $\{\hat{x},  \hat{y}, x,y,z\}^*$. Let $\rho$ be the
congruence on $\{x,y\}^*$ generated by $\{(u_i,v_i), \ i\in I\}$.
Let
\begin{enumerate}
\item[$(3)'$] $\lfloor zu\rfloor=\lfloor zv\rfloor,\ (u,v)\in \rho$ with $u>v$.
\end{enumerate}

\begin{theorem}\label{L4.12}With the above notation, the set $S_1=\{(1), (2), (3)'\}$ is a Gr\"{o}bner-Shirshov basis
in $Lie( \hat{x}, \hat{y}, x,y,z)$.
\end{theorem}
{\bf Proof:} For any $u\in \{x,y\}^*$, by induction on $|u|$,
$\overline{\lfloor zu\rfloor}=zu$. All possible compositions in
$S_1$ are intersection of (2) and $(3)'$, and inclusion of $(3)'$
and $(3)'$.

For $(2)\wedge(3)',\ w=\hat{x}zu,\ (u,v)\in \rho,\ u>v,\ f=[
\hat{x}z]+[zx],\ g=\lfloor zu\rfloor-\lfloor zv\rfloor$. We have
\begin{eqnarray*}
&&([ \hat{x}z]+[zx],\lfloor zu\rfloor-\lfloor
zv\rfloor)_w=[fu]_{\bar f}-[\hat{x}g]_{\bar g}\\
&\equiv&\lfloor ([ \hat{x}z]+[zx])u\rfloor
-[ \hat{x}(\lfloor zu\rfloor-\lfloor zv\rfloor)]\\
&\equiv &\lfloor zxu\rfloor+\lfloor  \hat{x}zv\rfloor\equiv\lfloor
zxu\rfloor-\lfloor zxv\rfloor\equiv0\ \ mod(S_1,w).
\end{eqnarray*}
For $(3)'\wedge(3)',\ w=zu_1=zu_2e,\ e\in \{x,y\}^*,\ (u_i,v_i)\in
\rho,\ u_i>v_i,\ i=1,2$. We have
\begin{eqnarray*}
&&(\lfloor zu_1\rfloor-\lfloor zv_1\rfloor,\lfloor
zu_2\rfloor-\lfloor zv_2\rfloor)_w\equiv
(\lfloor zu_1\rfloor-\lfloor zv_1\rfloor)-\lfloor(\lfloor zu_2\rfloor-\lfloor zv_2\rfloor)e\rfloor\\
&\equiv &\lfloor\lfloor zv_2\rfloor e\rfloor-\lfloor
zv_1\rfloor\equiv\lfloor zv_2e\rfloor-\lfloor zv_1\rfloor\equiv0\ \
mod(S_1,w).
\end{eqnarray*}

Thus, the set $S_1=\{(1), (2), (3)'\}$ is a Gr\"{o}bner-Shirshov
basis in $Lie( \hat{x}, \hat{y}, x,y,z)$. \hfill $\blacksquare$

\begin{corollary}(Kukin \cite{Kukin})\label{t4.12}
Let $u,v\in \{x,y\}^*$. Then
$$
u=v\ \mbox{ in the semigroup}\ P\Leftrightarrow\lfloor
zu\rfloor=\lfloor zv\rfloor\ \mbox{ in the Lie algebra }\ A_P.
$$
\end{corollary}
{\bf Proof:} Suppose that $ u=v\ \mbox{ in the semigroup}\ P$.
Without loss of generality, we may assume that $u=au_1b,\ v=av_1b$
for some $a,b\in \{x,y\}^*$ and $(u_1,v_1)\in \rho$. For any $r\in
\{x,y\}$, by the relations (1), we have $[\hat{x}r]=0$ and so
$\lfloor zxc\rfloor=\lfloor[z \hat{x}]c\rfloor=[\lfloor zc\rfloor
\hat{x}],\ \lfloor zyc\rfloor=[\lfloor zc\rfloor\hat{y}]$ for any
$c\in\{x,y\}^*$. From this it follows that in $A_P$, $\lfloor
zu\rfloor=\lfloor zau_1b\rfloor=\lfloor\lfloor zau_1\rfloor
b\rfloor=\lfloor\lfloor zu_1\widehat{\overleftarrow{a}}\rfloor
b\rfloor=\lfloor zu_1\widehat{\overleftarrow{a}}b\rfloor=\lfloor
zv_1\widehat{\overleftarrow{a}}b\rfloor=\lfloor
zav_1b\rfloor=\lfloor zv\rfloor$, where
$\overleftarrow{x_{i_1}x_{i_2}\cdots
x_{i_n}}=x_{i_n}x_{i_{n-1}}\cdots x_{i_1}$ and
$\widehat{x_{i_1}x_{i_2}\cdots
x_{i_n}}=\widehat{x_{i_1}}\widehat{x_{i_{2}}}\cdots
\widehat{x_{i_n}}$, $x_{i_j}\in \{x,y\}$. Moreover, $(3)'$ holds in
$A_P$.

Suppose that $\lfloor zu\rfloor=\lfloor zv\rfloor\ \mbox{ in the Lie
algebra }\ A_P$. Then both $\lfloor zu\rfloor$ and $\lfloor
zv\rfloor$  have the same normal form in $A_P$. Since $S_1$ is a
Gr\"{o}bner-Shirshov basis in $A_P$ by Theorem \ref{L4.12}, both
$\lfloor zu\rfloor$ and $\lfloor zv\rfloor$ can be reduced to the
same normal form of the form $\lfloor zc\rfloor$ for some $c\in
\{x,y\}^*$ only by the relations $(3)'$. This implies that in $P$,
$u=c=v$. \hfill $\blacksquare$

\ \

By the above corollary, if the semigroup $P$ has the undecidable
word  problem then so does the Lie algebra $A_P$.

\section{Gr\"{o}bner-Shirshov basis for the Drinfeld-Kohno Lie
algebra $\textbf{L}_n$}\label{dif3}

In this section we give a Gr\"{o}bner-Shirshov basis for the
Drinfeld-Kohno Lie algebra $\textbf{L}_n$.
\begin{definition}(\cite{[Et]})
Let $ n>2$ be an integer. The Drinfeld-Kohno Lie algebra
$\textbf{L}_{n}$ over $\mathbb{Z}$ is defined  by generators
$t_{ij}=t_{ji}$ for distinct indices $1 \leq i,\ j \leq n-1 $, and
relations
\begin{eqnarray*}
t_{ij}t_{kl}=0,  \\
t_{ij}(t_{ik}+t_{jk})= 0,
\end{eqnarray*}
where $i ,\  j , \ k ,\  l$ are distinct.
\end{definition}

Clearly, $\textbf{L}_{n}$ has a presentation $Lie_{\mathbb{Z}}(T|
S)$, where $T=\{t_{ij}|\ 1 \leq i<j\leq n-1\}$ and $S$ consists of
the following relations
\begin{eqnarray}\label{a}
&&t_{ij}t_{kl}=\ 0, \ \ \  \ k<i<j,\ k<l,\ l\neq\ i,\ j\\
&&\label{b} t_{jk}t_{ij}+t_{ik}t_{ij}=\ 0,\ \ \  \
i<j<k\\
&&\label{c} t_{jk}t_{ik}-t_{ik}t_{ij}=\ 0,\ \  \ i<j<k
\end{eqnarray}

Now we order $T$: $t_{ij}<t_{kl}$ if either $i<k$ or $i=k$ and
$j<l$. Let $<$ be the deg-lex order on $T^*$.

\begin{theorem}\label{t1} Let $S=\{(\ref{a}), (\ref{b}), (\ref{c})\}$ be as before, $<$ the deg-lex order on $T^*$. Then  $S$
is a Gr\"{o}bner-Shirshov basis for $\textbf{L}_{n}$.
\end{theorem}
\noindent\textbf{Proof.} We list all the possible ambiguities.
Denote  $(i)\wedge (j)$ the composition of the type $(i)$ and type
$(j)$.

For $(1)\wedge(n)$, $1\leq n\leq3$, the possible ambiguities $w$'s
are:
\begin{eqnarray*}
&(1)\wedge(1)& t_{ij}t_{kl}t_{mr}, \ (k<i<j,\ k<l,\  l\neq\
i,j,\ m<k<l,\ m<r, \ r\neq\ k,l),\\
&(1)\wedge(2)& t_{ij}t_{kl}t_{mk},\ (\ k<i<j,\
k<l,\ l\neq\ i,j,\ m<k<l),\\
&(1)\wedge(3)& t_{ij}t_{kl}t_{ml},\ ( k<i<j,\
k<l,\ l\neq\ i,j,\ m<k<l).\\
\end{eqnarray*}

For $(2)\wedge(n)$, $1\leq n\leq3$,
 the possible ambiguities $w$'s are:
\begin{eqnarray*}
&(2)\wedge(1)& t_{jk}t_{ij}t_{mr}, \ (\ m<i<j<k,\ m<r, \ r\neq\ i,j),\\
&(2)\wedge(2)& t_{jk}t_{ij}t_{mi},\ (\ m<i<j<k),\\
&(2)\wedge(3)& t_{jk}t_{ij}t_{mj},\ (\ m<i<j<k).\\
\end{eqnarray*}

For $(3)\wedge(n)$, $1\leq n\leq3$,
 the possible ambiguities $w$'s are:
 \begin{eqnarray*}
&(3)\wedge(1)& t_{jk}t_{ik}t_{mr}, \ (\ m<i<j<k,\ m<r, \ r\neq\ i,k),\\
&(3)\wedge(2)& t_{jk}t_{ik}t_{mi},\ (\ m<i<j<k),\\
&(3)\wedge(3)& t_{jk}t_{ik}t_{mk},\ (\ m<i<j<k).
\end{eqnarray*}

We claim that all compositions are trivial relative to $S$.

Here, we only prove cases $(1)\wedge(1)$, $(1)\wedge(2)$,
$(2)\wedge(1)$, $(2)\wedge(2)$, and the other cases can be proved
similarly.

 For $(1)\wedge(1)$,  let $f=t_{ij}t_{kl},\
g=t_{kl}t_{mr},\ k<i<j,\ k<l,\  l\neq\ i,j,\ m<k<l,\ m<r,\ r\neq\
k,l. $  Then $w=t_{ij}t_{kl}t_{mr}$  and
\begin{eqnarray*}
(f,g)_{w}&=&(t_{ij}t_{kl})t_{mr}-t_{ij}(t_{kl}t_{mr})\\
&=&\ (t_{ij}t_{mr})t_{kl}\  mod(S,w).
\end{eqnarray*}
There are three subcases to consider: $r\neq\ i,j$, $r=i$, $r=j$.

Subcase 1. If \ $r\neq\ i,j$,\ then\
\begin{eqnarray*}(t_{ij}t_{mr})t_{kl}\equiv\ 0\ mod(S,w).
\end{eqnarray*}

Subcase 2. If \ $r=i$,\ then
\begin{eqnarray*}
(t_{ij}t_{mr})t_{kl}&=&(t_{ij}t_{mi})t_{kl}\\
&\equiv&\ t_{kl}(t_{mj}t_{mi})\\  &\equiv&\
(t_{kl}t_{mj})t_{mi}+t_{mj}(t_{kl}t_{mi})\\ &\equiv&\  0\  mod(S,w).
\end{eqnarray*}

Subcase 3. If \ $r=j$,\ then
\begin{eqnarray*}
(t_{ij}t_{mr})t_{kl}&=&(t_{ij}t_{mj})t_{kl}\\
&\equiv&\ -t_{kl}(t_{mj}t_{mi})\\  &\equiv&\
-(t_{kl}t_{mj})t_{mi}-t_{mj}(t_{kl}t_{mi})\\ &\equiv&\  0\ mod(S,w).
\end{eqnarray*}

For $(1)\wedge(2)$,  let $f=t_{ij}t_{kl},\
g=t_{kl}t_{mk}+t_{ml}t_{mk},\ k<i<j,\ k<l,\ l\neq\ i,j,\ m<k<l.\
$Then $w=t_{ij}t_{kl}t_{mk}$ and
\begin{eqnarray*}
(f,g)_{w}&=&(t_{ij}t_{kl})t_{mk}-t_{ij}(t_{kl}t_{mk}+t_{ml}t_{mk})\\
&=&\ (t_{ij}t_{mk})t_{kl}-t_{ij}(t_{ml}t_{mk})\\ &\equiv&\
-t_{ij}(t_{ml}t_{mk})\\ &\equiv&\
-(t_{ij}t_{ml})t_{mk}-t_{ml}(t_{ij}t_{mk})\\ &\equiv&\  0\ mod(S,w).
\end{eqnarray*}

For  $(2)\wedge(1)$,  let $f=t_{jk}t_{ij}+t_{ik}t_{ij},\
g=t_{ij}t_{mr},\ m<i<j<k,\ \ m<r,\ r\neq\ i,j.\  $Then
$w=t_{jk}t_{ij}t_{mr}$ and
\begin{eqnarray*}
(f,g)_{w}&=&(t_{jk}t_{ij}+t_{ik}t_{ij})t_{mr}-t_{jk}(t_{ij}t_{mr})\\
&\equiv&\ (t_{jk}t_{mr})t_{ij}+(t_{ik}t_{mr})t_{ij}\  mod(S,w).
\end{eqnarray*}
There are two subcases to consider: $r\neq\ k$, $r=k$. \\

Subcase 1. If \ $r\neq\ k$,\ then\
\begin{eqnarray*}(t_{jk}t_{mr})t_{ij}+(t_{ik}t_{mr})t_{ij}\equiv\ 0\ mod(S,w).
\end{eqnarray*}

Subcase 2. If \ $r=k$,\ then
\begin{eqnarray*}
(t_{jk}t_{mr})t_{ij}+(t_{ik}t_{mr})t_{ij}&=&(t_{jk}t_{mk})t_{ij}+(t_{ik}t_{mk})t_{ij}\\
&\equiv&\ -t_{ij}(t_{mk}t_{mj})-t_{ij}(t_{mk}t_{mi})\\  &\equiv&\
(t_{ij}t_{mj})t_{mk}+(t_{ij}t_{mi})t_{mk}\\ &\equiv&\
-t_{mk}(t_{mj}t_{mi})+t_{mk}(t_{mj}t_{mi}) \\ &\equiv&\ 0\ mod(S,w).
\end{eqnarray*}

For $(2)\wedge(2)$,  let $f=t_{jk}t_{ij}+t_{ik}t_{ij},\
g=t_{ij}t_{mi}+t_{mj}t_{mi},\ m<i<j<k.\ $Then
$w=t_{jk}t_{ij}t_{mi}$\ and
\begin{eqnarray*}
(f,g)_{w}&=&(t_{jk}t_{ij}+t_{ik}t_{ij})t_{mi}-t_{jk}(t_{ij}t_{mi}+t_{mj}t_{mi})\\
&=&\
(t_{jk}t_{mi})t_{ij}+(t_{ik}t_{mi})t_{ij}+t_{ik}(t_{ij}t_{mi})-t_{jk}(t_{mj}t_{mi})\\
&\equiv&\
t_{ij}(t_{mk}t_{mi})-t_{ik}(t_{mj}t_{mi})-t_{jk}(t_{mj}t_{mi})\\
&\equiv&\
-(t_{ij}t_{mi})t_{mk}+(t_{ik}t_{mi})t_{mj}-(t_{jk}t_{mj})t_{mi}\\
&\equiv&\
-t_{mk}(t_{mj}t_{mi})-(t_{mk}t_{mi})t_{mj}+(t_{mk}t_{mj})t_{mi}\\
&\equiv&\  0\ mod(S,w).
\end{eqnarray*}

So $S$  is a Gr\"{o}bner-Shirshov basis for $\textbf{L}_{n}$. \hfill
$\blacksquare$

\ \

Let $L$ be a Lie algebra over a commutative ring $K$, $L_{1}$ an
ideal of $L$ and $L_{2}$ a subalgebra of $L$. We call $L$  a
semidirect product of $L_{1}$ and $L_{2}$ if $L=L_{1}\oplus L_{2}$
as $K$-modules.

By  Theorems \ref{cdLie} and \ref{t1}, we have immediately the
following corollaries.
\begin{corollary}\label{cor1}
The Drinfeld-Kohno Lie algebra $\textbf{L}_{n}$ is a free
$\mathbb{Z}$-module with a $\mathbb{Z}$-basis
$$Irr(S)=\{[t_{ik_{1}}t_{ik_{2}}\cdots\ t_{ik_{m}}]\ |\ t_{ik_{1}}t_{ik_{2}}\cdots\ t_{ik_{m}}\
is \ an \mbox { ALSW } in \ T^*,\ m\in\mathbb{N}\}.$$
\end{corollary}

\begin{corollary}\label{cor1}(\cite{[Et]})
$\textbf{L}_{n}$ is an iterated semidirect product of free Lie
algebras.
\end{corollary}

\noindent\textbf{Proof.} \ Let $A_{i}$ be the free Lie algebra
generated by $\{t_{ij}\ |\ i<j\leq\ n-1\}$. Clearly,
$$\textbf{L}_{n}=A_{1}\oplus\ A_{2}\oplus\ \cdots\ \oplus\ A_{n-2}$$
as $\mathbb{Z}$-modules, and from the relations $(1),(2),(3)$, we
have
$$
A_{i}\triangleleft\ A_{i}+A_{i+1}+\ \cdots\ +A_{n-2}.
$$
\hfill$\blacksquare$

\ \

{\bf Remark:} In this section, if we replace the base ring
$\mathbb{Z}$ by an arbitrary commutative ring $K$ with identity,
then all results hold.


\begin{thebibliography}{99}


\bibitem{AL}W.W. Adams, P. Loustaunau, An introduction to Gr\"{o}bner
bases, Graduate Studies in Mathematics, Vol. 3, American
Mathematical Society (AMS), 1994.


\bibitem{Be78}G.M. Bergman, The diamond lemma for ring theory, {\it Adv.
Math.} 29 (1978) 178-218.


\bibitem{Bo72} L.A. Bokut, Insolvability of the word problem for Lie algebras and
subalgebras of finitely presented Lie algebras,  Izvestija AN USSR
(mathem.) 36 (1972) 1173-1219.


\bibitem{Bo76}L.A. Bokut, Imbeddings into simple associative
algebras, {\it Algebra i Logika} 15 (1976) 117-142.



\bibitem{BC}L.A. Bokut, Y.Q. Chen, Gr\"{o}bner-Shirshov
 bases: Some new results, Proceedings of the Second International Congress
 in Algebra and Combinatorics, World Scientific, (2008) 35-56.

\bibitem{BC13}L.A. Bokut, Y.Q. Chen, Gr\"obner-Shirshov bases and their calculation,
arxiv.org/abs/1303.5366.


\bibitem{BCC08}L.A. Bokut, Y.Q. Chen, Y.S. Chen, Composition-Diamond lemma
for tensor product of free algebras, \emph{Journal of Algebra} 323
(2010) 2520-2537.

\bibitem{BCC11}L.A. Bokut, Y.Q. Chen, Y.S. Chen,  Groebner-Shirshov bases for Lie
algebras over a commutative algebra,  \emph{Journal of Algebra} 337
(2011) 82-102.


\bibitem{BCD08}L.A. Bokut, Y.Q. Chen, X.M. Deng,
Gr\"{o}bner-Shirshov bases for Rota-Baxter algebras, \emph{Siberian
Math. J.} 51 (6) (2010) 978-988.

\bibitem{BCLi13}L.A. Bokut,  Y.Q. Chen, Y. Li,
Lyndon-Shirshov basis and anti-commutative algebras, \emph{Journal
of Algebra} 378 (2013) 173-183.


\bibitem{BCL08}L.A. Bokut,  Y.Q. Chen, C.H. Liu,
Gr\"{o}bner-Shirshov bases for dialgebras,\emph{ International
Journal of Algebra and Computation} 20 (3) (2010) 391-415.

\bibitem{BCM} L.A. Bokut, Y.Q. Chen, Q.H. Mo,
Gr\"{o}bner-Shirshov bases and embeddings of algebras,
\emph{International Journal of Algebra and Computation} 20 (2010)
875-900.


\bibitem{BCM13} L.A. Bokut, Y.Q. Chen, Q.H. Mo, Gr\"{o}bner-Shirshov bases for
semirings, \emph{Journal of Algebra} 385 (2013) 47-63.


\bibitem{BCS} L.A. Bokut, Y.Q. Chen,
K.P. Shum, Some new results on Groebner-Shirshov bases, in:
Proceedings of International Conference on Algebra 2010, Advances in
Algebraic Structures, (2012) 53-102.



\bibitem{BFKK00}L.A. Bokut, Y. Fong, W.-F. Ke, P.S. Kolesnikov, Gr\"obner and
Gr\"obner-Shirshov bases in algebra and conformal algebras, {\it
Fundamental and Applied Mathematics} 6 (3) (2000) 669-706.

\bibitem{BK03}L.A. Bokut, P.S. Kolesnikov, Gr\"obner-Shirshov bases: from their
incipiency to the present, {\it J. Math. Sci.} 116 (1) (2003)
2894-2916.


\bibitem{BK05}L.A. Bokut, P.S. Kolesnikov, Gr\"obner-Shirshov bases, conformal
algebras and pseudo-algebras, {\it J. Math. Sci.} 131 (5) (2005)
5962-6003.


\bibitem{BKu94}L.A. Bokut, G. Kukin, Algorithmic and Combinatorial algebra, Kluwer Academic Publ., Dordrecht,
1994


\bibitem{Boon} W.W. Boone, The word problem, {\it Ann. Math}. 70
(1959) 207-265.

\bibitem{Bu70}B. Buchberger, An algorithmical criteria for the
solvability of algebraic systems of equations, {\it Aequationes
Math.} 4 (1970) 374-383.

\bibitem{BuCL}B. Buchberger, G.E. Collins, R. Loos, R. Albrecht, Computer algebra,
symbolic and algebraic computation, Computing Supplementum, Vol.4,
New York: Springer-Verlag, 1982.

\bibitem{BuW}B. Buchberger, Franz Winkler, Gr\"{o}bner bases and
applications, London Mathematical Society Lecture Note Series,
Vol.251, Cambridge: Cambridge University Press, 1998.




\bibitem{CC12}Y.S. Chen, Y.Q. Chen,  Groebner-Shirshov bases for matabelian Lie
algebras,  \emph{Journal of Algebra} 358 (2012) 143-161.




\bibitem{CLO} D.A. Cox, J. Little, D. O'Shea, Ideals,
varieties and algorithms: An introduction to computational algebraic
geometry and commutative algebra, Undergraduate Texts in
Mathematics, New York: Springer-Verlag, 1992.


\bibitem{Ei} D. Eisenbud, Commutative algebra with a view toward algebraic geometry,
Graduate Texts in Math., Vol.150, Berlin and New York:
Springer-Verlag, 1995.

\bibitem{[Dr]}V.G. Drinfeld, Quasi-Hopf algebras, {\it Algebra i Analiz}
1 (1989) 114-148.

\bibitem{[Et]}P. Etingof, A. Henriques, J. Kamnitzer, E.M. Rains,
The cohomology ring of the real locus of the moduli space of stable
curves of genus 0 with marked points, {\it Ann. Math.} 171 (2010)
731-777.


\bibitem{Higman} G. Higman,  Subgroups of finitely presented groups.
{\it Proc. Royal Soc. London (Series A)} 262 (1961) 455-475.




\bibitem{Hi64}H. Hironaka, Resolution of singularities of an algebraic variety
over a field if characteristic zero, I, II, {\it Ann.  Math.} 79
(1964) 109-203, 205-326.


\bibitem{KL}S.-J. Kang, K.-H. Lee, Gr\"obner-Shirshov bases for irreducible
$sl_{n+1}$-modules, {\it Journal of Algebra} 232 (2000) 1-20.


\bibitem{[Ko]}T. Kohno, Serie de Poincare Koszul associee aux groupes de tresses pures,
{\it Invent. Math.} 82 (1985) 57-75.


\bibitem{Kukin} G.P. Kukin, On the word problem for Lie algebras,
{\it Sibirsk. Mat. Zh.} 18 (1977) 1194-1197.

\bibitem{Lyndon}  Lyndon, R.C.: On Burnside's problem I, {\it  Trans. Amer. Math. Soc.} 77
(1954) 202-215.

\bibitem{Markov}  A.A. Markov, Impossibility of some algorithms in the
theory of some associative system, {\it Dokl. Akad. Nauk SSSR} 55
(1947) 587-590.


\bibitem{MZ}A.A. Mikhalev, A.A. Zolotykh, Standard Gr\"obner-Shirshov bases of free algebras over
rings, I. Free associative algebras,  \emph{International Journal of
Algebra and Computation} 8 (6) (1998) 689-726.



\bibitem{Novikov} P.S. Novikov, On algorithmic undecidability of the
word problem in the theory of groups,  Trudy MKat. Inst. Steklov. 44
(1955) 1-144.

\bibitem{Post46}E. Post, A variant of a recursively unsolvable
problem, {\it Bull. Amer. Math. Soc.} 52 (1946) 264-268.



\bibitem{Sh58} A.I. Shirshov,   On free Lie rings,  {\it Mat. Sb.} 45 (2) (1958) 113-122.




\bibitem{Sh62b} A.I. Shirshov, Some algorithmic problem for Lie
algebras, {\it Sibirsk. Mat. Zh.} 3 (2) (1962) 292-296; English
translation in SIGSAM Bull.  33 (1999) 3-6.


\bibitem{Sh62a}A.I. Shirshov, Some algorithmic problem for $\varepsilon$-algebras,
 {\it Sibirsk. Mat. Z.} 3 (1962) 132-137.

\bibitem{Shir3} Selected works of A.I. Shirshov, Eds L.A. Bokut, V. Latyshev, I. Shestakov,
E. Zelmanov, Trs M. Bremner, M. Kochetov, Birkh\"auser, Basel,
Boston, Berlin,  2009.

\bibitem{Turing} A.M. Turing, The word problem in semi-groups with
cancellation,  {\it Ann.  Math.} 52 (1950) 191-505.

\end{thebibliography}
\end{document}